\documentclass[12pt,thmsb]{article}
\usepackage{amsfonts}
\usepackage{amsmath}
\usepackage{amssymb}
\usepackage{sw20jart}

\input{tcilatex}
\input tcilatex
\begin{document}

\author{Refik Keskin$^{a}$ and Olcay Karaatl\i $^{a}$ \and %
rkeskin@sakarya.edu.tr, okaraatli@sakarya.edu.tr \\
$^{a}$Sakarya University, Faculty of Arts and Science, \\
Department of Mathematics, Sakarya, Turkey}
\title{Generalized Fibonacci and Lucas Numbers of the form $5x^{2}$}
\maketitle

\begin{abstract}
Let $\left( U_{n}(P,Q\right) $ and $\left( V_{n}(P,Q\right) $ denote the
generalized Fibonacci and Lucas sequence, respectively. In this study, we
assume that $Q=1.$ We determine all indices $n$ such that $U_{n}=5\square $
and $U_{n}=5U_{m}\square $ under some assumptions on $P.$ We show that the
equation $V_{n}=5\square $ has the solution only if $n=1$ for the case when $%
P$ is odd. Moreover, we show that the equation $V_{n}=5V_{m}\square $ has no
solutions.
\end{abstract}


\emph{MSC: }$11$B$37,$ $11$B$39,$ $11$A$07$

\section{Introduction}

Let $P$ and $Q$ be non-zero integers with $P^{2}+4Q\neq 0.$ Generalized
Fibonacci sequence $\left( U_{n}(P,Q)\right) $ and Lucas sequence $\left(
V_{n}(P,Q)\right) $ are defined by the following recurrence relations:%
\begin{equation*}
U_{0}(P,Q)=0,\text{ }U_{1}(P,Q)=1,\text{ }%
U_{n+2}(P,Q)=PU_{n+1}(P,Q)+QU_{n}(P,Q)\text{ for }n\geq 0
\end{equation*}%
and 
\begin{equation*}
V_{0}(P,Q)=2,\text{ }V_{1}(P,Q)=P,\text{ }%
V_{n+2}(P,Q)=PV_{n+1}(P,Q)+QV_{n}(P,Q)\text{ for }n\geq 0.
\end{equation*}%
$U_{n}(P,Q)$ is called the $n-$th generalized Fibonacci number and $%
V_{n}(P,Q)$ is called the $n-$th generalized Lucas number. Also generalized
Fibonacci and Lucas numbers for negative subscripts are defined as 
\begin{equation*}
U_{-n}(P,Q)=\frac{-U_{n}(P,Q)}{(-Q)^{n}}\text{ and }V_{-n}=\frac{V_{n}(P,Q)}{%
(-Q)^{n}}\text{ for }n\geq 1,
\end{equation*}%
respectively. Taking $\alpha =(P+\sqrt{P^{2}+4Q})/2$ and $\beta =(P-\sqrt{%
P^{2}+4Q})/2$ to be the roots of the characteristic equation $x^{2}-Px-Q=0,$
we have the well known expressions named Binet forms%
\begin{equation*}
U_{n}(P,Q)=(\alpha ^{n}-\beta ^{n})/(\alpha -\beta )\text{ and }%
V(P,Q)=\alpha ^{n}+\beta ^{n}
\end{equation*}%
for all $n\in 
\mathbb{Z}
.$ From now on, we assume that $P>0$ and $P^{2}+4Q$ $>0$. Instead of $%
U_{n}(P,Q)$ and $V_{n}(P,Q),$ we will use $U_{n}$ and $V_{n},$ respectively.
For $P=Q=1,$ the sequence $\left( U_{n}\right) $ is the familiar Fibonacci
sequence $\left( F_{n}\right) $ and the sequence $\left( V_{n}\right) $ the
familiar Lucas sequence $\left( L_{n}\right) .$ If $P=2$ and $Q=1,$ then we
have the well known Pell sequence $\left( P_{n}\right) $ and Pell-Lucas
sequence $\left( Q_{n}\right) .$ For $Q=-1,$ we represent $\left(
U_{n}\right) $ and $\left( V_{n}\right) $ by $\left( u_{n}\right) $ and $%
\left( v_{n}\right) ,$ respectively. For more information about generalized
Fibonacci and Lucas sequences, one can consult \cite{KLMN, MSKT, RIBENBO,
RABINO}.

Investigations of the properties of second order linear recurring sequences,
have given rise to questions concerning whether, for certain pairs $(P,Q),$ $%
U_{n}$ or $V_{n}$ is square(=$\square $). In particular, the squares in
sequences $\left( U_{n}\right) $ and $\left( V_{n}\right) $ were
investigated by many authors. Ljunggrenn \cite{LJUN} showed in 1942 that if $%
P=2,$ $Q=1,$ and $n\geq 2,$ then $U_{n}=\square $ precisely for $n=7$ and $%
U_{n}=2\square $ precisely for $n=2.$ In 1964, Cohn \cite{CO1} proved that
if $P=Q=1,$ then the only perfect square greater than $1$ in the sequence $%
\left( U_{n}\right) $ is $U_{12}=12^{2}$ $($see also Alfred \cite{ALF}, Burr 
\cite{BURR}, and Wyler \cite{WYLER}), and he \cite{CO2, CO3} solved the
equations $U_{n}=2\square $ and $V_{n}=\square ,2\square .$ Furthermore, in
other papers, Cohn \cite{CO4, CO5} determined the squares and twice the
squares in $\left( U_{n}\right) $ and $\left( V_{n}\right) $ when $P$ is odd
and $Q=\pm 1.$ Ribenboim and McDaniel \cite{MC DAN1} determined all indices $%
n$ such that $U_{n}=\square ,$ $2U_{n}=\square ,$ $V_{n}=\square $ or $%
2V_{n}=\square $ for all odd relatively prime integers $P$ and $Q.$ In 1998,
Kagawa and Terai \cite{KGW} considered a similar problem for the case when $%
P $ is even and $Q=1.$ Using the elementary properties of elliptic curves,
they showed that if $P=2t$ with $t$ even and $Q=1,$ then $U_{n}=\square ,$ $%
2U_{n}=\square ,$ $V_{n}=\square $ or $2V_{n}=\square $ implies $n\leq 3$
under some assumptions. Besides, for $Q=1,$ Nakamula and Petho \cite{PETHO}
gave the solutions of the equations $U_{n}=w\square $ where $w\in \left\{
1,2,3,6\right\} .$ In 1998, Ribenboim and McDaniel \cite{MC DAN2} showed
that if $P$ is even, $Q\equiv 3(\func{mod}$ $4)$ and $U_{n}=\square ,$ then $%
n$ is a square or twice an odd square and all prime factors of $n$ divides $%
P^{2}+4Q.$ In a latter paper, the same authors \cite{MC DAN3} solved the
equation $V_{n}=3\square $ for $P\equiv 1,3(\func{mod}$ $8),Q\equiv 3(\func{%
mod}$ $4),$ $(P,Q)=1$ and solved the equation $U_{n}=3\square $ for all odd
relatively prime integers $P$ and $Q.$ Moreover, in \cite{CO6}, Cohn solved
the equations $V_{n}=V_{m}x^{2}$ and $V_{n}=2V_{m}x^{2}$ when $P$ is odd.
Keskin and Yosma \cite{KSKN} gave the solutions of the equations $%
F_{n}=2F_{m}\square ,$ $L_{n}=2L_{m}\square ,$ $F_{n}=3F_{m}\square ,$ $%
F_{n}=6F_{m}\square ,$ and $L_{n}=6L_{m}\square .$ In \cite{ZAFER}, \c{S}iar
and Keskin, assuming $Q=1,$ solve the equation $V_{n}=2V_{m}\square $ when $%
P $ is even. They determine all indices $n$ such that $V_{n}=kx^{2}$ when $%
k|P$ and $P$ is odd, where $k$ is a square-free positive divisor of $P.$
They show that there is no integer solution of the equations $V_{n}=3\square 
$ and $V_{n}=6\square $ for the case when $P$ is odd and also they give the
solution of the equations $V_{n}=3V_{m}\square $ and $V_{n}=6V_{m}\square .$
More generally, we can give the following theorem proved by Shorey and
Stewart in \cite{SHOREY}:

Let $\ A>0$ be an integer. Then there exists an effectively computable
number $C>0,$ which depends on $A,$ such that if $n>0$ and $U_{n}=A\square $
or $V_{n}=A\square ,$ then $n<C.$

In this study, we assume, from this point on, that $Q=1.$ We determine all
indices $n$ such that $U_{n}=5\square $ and $U_{n}=5U_{m}\square $ under
some assumptions on $P.$ We show that if $P$ is odd, then the equation $%
V_{n}=5\square $ has the solution only if $n=1.$ Moreover, we prove that the
equation $V_{n}=5V_{m}\square $ has no solutions.

\section{Preliminaries}

In this section, we give some theorems, lemmas and well known identities
about generalized Fibonacci and Lucas numbers, which will be needed in the
proofs of the main theorems. Through the paper $(\frac{\ast }{\ast })$
denotes the Jacobi symbol. The proofs of the following two theorems can be
found in \cite{KSKN1}.

\begin{theorem}
\label{t2.1}Let $m,r\in \mathbb{%
\mathbb{Z}
}$ and $n$ be non-zero integer. Then
\end{theorem}

\begin{equation}
U_{2mn+r}\equiv \left( -1\right) ^{mn}U_{r}\left( \func{mod}\text{ }%
U_{m}\right)  \label{2.1}
\end{equation}%
and%
\begin{equation}
V_{2mn+r}\equiv \left( -1\right) ^{mn}V_{r}\left( \func{mod}\text{ }%
U_{m}\right)  \label{2.2}
\end{equation}

\begin{theorem}
\label{t2.2}Let $m,r\in \mathbb{%
\mathbb{Z}
}$ and $n$ be non-zero integer. Then
\end{theorem}

\begin{equation}
U_{2mn+r}\equiv \left( -1\right) ^{(m+1)n}U_{r}\left( \func{mod}\text{ }%
V_{m}\right)  \label{2.3}
\end{equation}%
and%
\begin{equation}
V_{2mn+r}\equiv \left( -1\right) ^{(m+1)n}V_{r}\left( \func{mod}\text{ }%
V_{m}\right)  \label{2.4}
\end{equation}%
We state the following theorem from \cite{PETHO}.

\begin{theorem}
\label{t2.3}Let $P>0$ and $Q=1.$ If $U_{n}=wx^{2}$ with $w\in \left\{
1,2,3,6\right\} ,$ then $n\leq 2$ except when $%
(P,n,w)=(2,4,3),(2,7,1),(4,4,2),(1,12,1),(1,3,2),(1,4,3),(1,6,2),$ and $%
(24,4,3).$
\end{theorem}

We give the following two theorems from \cite{CO4} and \cite{CO5}.

\begin{theorem}
\label{t2.4}If $P$ is odd, then the equation $V_{n}=x^{2}$ has the solutions 
$n=1,$ $P=\square ,$ and $P\neq 1$ or $n=1,3$ and $P=1$ or $n=3$ and $P=3.$
\end{theorem}

\begin{theorem}
\label{t2.5}If $P$ is odd, then the equation $V_{n}=2x^{2}$ has the
solutions $n=0$ or $n=\pm 6$ and $P=1,5.$
\end{theorem}

The following two theorems can be obtained from Theorem $11$ and Theorem $12$
given in \cite{CO6}.

\begin{theorem}
\label{T2.6}Let $P$ be an odd integer, $m\geq 1$ be an integer and $%
V_{n}=V_{m}x^{2}$ for some integer $x.$ Then $n=m.$
\end{theorem}

\begin{theorem}
\label{T2.7}If $P$ is an odd integer and $m\geq 1,$ then there is no integer 
$x$ such that $V_{n}=2V_{m}x^{2}.$
\end{theorem}

The following theorem can be obtained from Theorem $6$ given in \cite{CO6}.

\begin{theorem}
\label{T2.9}Let $P$ be an odd integer, $m\geq 2$ be an integer and $%
U_{n}=2U_{m}x^{2}$ for some integer $x.$ Then $n=12,m=6,P=5.$
\end{theorem}

Now we give some well known theorems in number theory. For more detailed
information, see \cite{NIV} or \cite{BURT}.

\begin{theorem}
\label{t2.6}Let $m$ be an odd integer. Suppose that $x^{2}\equiv -a^{2}(%
\func{mod}$ $m)$ for some nonzero integers $x$ and $a.$ Then $m\equiv 1(%
\func{mod}$ $4).$
\end{theorem}

We omit the proof of the following theorem since it can be seen easily by
induction method.

\begin{theorem}
\label{t2.7}Let $k$ be an integer with $k\geq 1.$ Then $L_{2^{k}}\equiv 3(%
\func{mod}$ $4).$
\end{theorem}

\begin{corollary}
\label{c2.1}Let $a$ be any nonzero integer. If $k\geq 1,$ then there is no
integer $x$ such that $x^{2}\equiv -a^{2}(\func{mod}$ $L_{2^{k}}).$
\end{corollary}

We omit the proof of the following theorem due to Keskin and Demirt\"{u}rk 
\cite{KSKN2}.

\begin{theorem}
\label{t2.8}All nonnegative integer solutions of the equation $%
u^{2}-5v^{2}=1 $ are given by $(u,v)=(L_{3z}/2,F_{3z}/2)$ with nonnegative
even integer $z$ and all nonnegative integer solutions of the equation $%
u^{2}-5v^{2}=-1$ are given by $(u,v)=(L_{3z}/2,F_{3z}/2)$ with positive odd
integer $z.$
\end{theorem}

By using the above theorem, we can give the following theorem without proof.

\begin{theorem}
\label{t1}All nonnegative integer solutions of the equation $%
x^{2}-4xy-y^{2}=-5$ are given by $(x,y)=(L_{3z+3}/2,L_{3z}/2)$ with
nonnegative even integer $z$ and all nonnegative integer solutions of the
equation $x^{2}-4xy-y^{2}=-1$ are given by $(x,y)=(F_{3z+3}/2,F_{3z}/2)$
with positive odd integer $z.$
\end{theorem}

The following lemma can be found in \cite{MC DAN3}.

\begin{lemma}
\label{l2.2}Let $P$ be odd, $m$ be an odd positive integer, and $r\geq 1.$
Then 
\begin{equation*}
V_{2^{r}m}\equiv \left\{ 
\begin{array}{c}
2\text{ }(\func{mod}\text{ }8)\text{ if }3\mid m, \\ 
3\text{ }(\func{mod}\text{ }8)\text{ if }3\nmid m\text{ and }r=1, \\ 
7\text{ }(\func{mod}\text{ }8)\text{ if }3\nmid m\text{ and }r>1.%
\end{array}%
\right.
\end{equation*}
\end{lemma}

Now we give the following results involving Fibonacci and Lucas numbers with
nonnegative integers $a$ and $m.$%
\begin{equation}
F_{m}=a^{2}\text{ iff }m=0,1,2,12,  \label{2.5}
\end{equation}%
\begin{equation}
F_{m}=2a^{2}\text{ iff }m=0,3,6,  \label{2.6}
\end{equation}%
\begin{equation}
F_{m}=5a^{2}\text{ iff }m=0,5,  \label{2.7}
\end{equation}%
\begin{equation}
F_{m}=10a^{2}\text{ iff }m=0,  \label{2.8}
\end{equation}%
\begin{equation}
L_{m}=a^{2}\text{ iff }m=1,3,  \label{2.9}
\end{equation}%
\begin{equation}
L_{m}=2a^{2}\text{ iff }m=0,6.  \label{2.10}
\end{equation}%
The equations (\ref{2.5}) and (\ref{2.6}) are Theorems $3$ and $4$ in \cite%
{CO3}; (\ref{2.7}) follows from Theorem $3$ in \cite{ROBBINS}; (\ref{2.8})
follows from Theorem $3$ in \cite{ROBBINS1}; (\ref{2.9}) and (\ref{2.10})
are Theorems $1$ and $2$ in \cite{CO3}.\qquad

We will need the following identities concerning generalized Fibonacci and
Lucas numbers:%
\begin{equation}
U_{2n}=U_{n}V_{n},  \label{2.11}
\end{equation}%
\begin{equation}
V_{2n}=V_{n}^{2}-2(-1)^{n},  \label{2.12}
\end{equation}%
\begin{equation}
V_{n}^{2}-(P^{2}+4)U_{n}^{2}=4(-1)^{n},  \label{2.13}
\end{equation}%
\begin{equation}
U_{3n}=U_{n}\left( (P^{2}+4)U_{n}^{2}+3(-1)^{n}\right) ,  \label{2.14}
\end{equation}%
\begin{equation}
u_{3n}=u_{n}\left( (P^{2}-4)u_{n}^{2}+3\right)  \label{10}
\end{equation}%
\begin{equation}
U_{5n}=\left\{ 
\begin{array}{c}
U_{n}\left( (P^{2}+4)^{2}U_{n}^{4}+5(P^{2}+4)U_{n}^{2}+5\right) \text{ if }n%
\text{ is even} \\ 
U_{n}\left( (P^{2}+4)^{2}U_{n}^{4}-5(P^{2}+4)U_{n}^{2}+5\right) \text{ if }n%
\text{ is odd,}%
\end{array}%
\right.  \label{2.15}
\end{equation}%
\begin{equation}
V_{5n}=\left\{ 
\begin{array}{c}
V_{n}(V_{n}^{4}-5V_{n}^{2}+5)\text{ if }n\text{ is even} \\ 
V_{n}(V_{n}^{4}+5V_{n}^{2}+5)\text{ if }n\text{ is odd,}%
\end{array}%
\right.  \label{2.16}
\end{equation}%
\begin{equation}
\text{If }m\geq 1,\text{ then }V_{m}|V_{n}\text{ iff }m|n\text{ and }n/m%
\text{ is odd integer,}  \label{2.17}
\end{equation}%
\begin{equation}
\text{If }U_{m}\neq 1,\text{ then }U_{m}|U_{n}\text{ iff }m|n.  \label{2.01}
\end{equation}%
\begin{equation}
\text{If }P\text{ is odd, then }(U_{n},V_{n})=\left\{ 
\begin{array}{c}
1\text{ if }3\nmid n \\ 
2\text{ if }3\mid n,%
\end{array}%
\right.  \label{2.18}
\end{equation}%
\begin{equation}
\text{If }r\geq 3,\text{ then }V_{2^{r}}\equiv 2(\func{mod}\text{ }V_{2}).
\label{2.22}
\end{equation}%
If $5|P$ and $n$ is odd, then $5|V_{n}$ and therefore from (\ref{2.16}), it
follows that $\ $%
\begin{equation}
V_{5n}=5V_{n}(5a+1)  \label{2.101}
\end{equation}%
for some positive integer $a.$ \qquad \qquad \qquad

\section{Main Theorems}

From this point on, we assume that $m,n\geq 1.$ Now we prove two theorems
which help us to determine for what values of $n,$ the equation $U_{n}=5x^{2}
$ has solutions and for what values of $m,n,$ the equations $%
V_{n}=5V_{m}x^{2}$ and $U_{n}=5U_{m}x^{2}$ have solutions.

\begin{theorem}
\label{t3.1}The only positive integer solution of the equation $%
x^{4}+3x^{2}+1=5y^{2}$ is given by $(x,y)=(1,1)$ and the only positive
integer solution of the equation $x^{4}-3x^{2}+1=5y^{2}$ is given by $%
(x,y)=(2,1).$
\end{theorem}

\proof%
Assume that $x^{4}\pm 3x^{2}+1=5y^{2}$ for some positive integers $x$ and $%
y. $ Multiplying both sides of the equations by $4$ and completing the
square gives 
\begin{equation*}
(2x\pm 3)^{2}-5=5(2y)^{2}.
\end{equation*}%
Then it follows that%
\begin{equation*}
(2y)^{2}-5\left( (2x\pm 3)/5\right) ^{2}=-1.
\end{equation*}%
By Theorem \ref{t3.1}, we get $2y=L_{3z}/2$ and $(2x^{2}\pm 3)/5=F_{3z}/2$
with positive odd integer $z.$ Assume that $z>1.$ Then we can write $z=4q\pm
1$ for some $q>0$ and therefore $z=2.2^{k}a\pm 1$ with $2\nmid a$ and $k\geq
1.$ Thus by (\ref{2.3}), we get 
\begin{equation*}
F_{3z}=F_{3(4q\pm 1)}=F_{12q\pm 3}=F_{2.2^{k}3a\pm 3}\equiv -F_{\pm 3}\equiv
-F_{3}(\func{mod}\text{ }L_{2^{k}}),
\end{equation*}%
i.e.,%
\begin{equation*}
F_{3z}\equiv -2(\func{mod}\text{ }L_{2^{k}}).
\end{equation*}%
Substituting the value of $F_{3z}$ and rewriting the above congruence gives 
\begin{equation*}
4x^{2}\pm 6\equiv -10(\func{mod}\text{ }L_{2^{k}}).
\end{equation*}%
This shows that 
\begin{equation*}
4x^{2}+6\equiv -10(\func{mod}\text{ }L_{2^{k}})\text{ or }4x^{2}-6\equiv -10(%
\func{mod}\text{ }L_{2^{k}}).
\end{equation*}%
Then it follows that 
\begin{equation*}
x^{2}\equiv -4(\func{mod}L_{2^{k}})
\end{equation*}%
or 
\begin{equation*}
x^{2}\equiv -1(\func{mod}L_{2^{k}}),
\end{equation*}%
which is a contradiction by Corollary \ref{c2.1}. Thus $z=1$ and therefore $%
2x^{2}\pm 3=5F_{3}/2$ and $2y=L_{3}/2.$ A simple computation shows that $y=1$
and $x=1$ or $x=2.$ This means that the equation $x^{4}+3x^{2}+1=5y^{2}$ has
only the positive integer solution $(x,y)=(1,1)$ and the equation $%
x^{4}-3x^{2}+1=5y^{2}$ has only the positive integer solution $(x,y)=(2,1).$%
This completes the proof of Theorem \ref{t3.1}.%
\endproof%

\begin{theorem}
\label{t3.2}The equation $x^{4}+5x^{2}+5=5y^{2}$ has no solution in positive
integers $x$ and $y.$
\end{theorem}

\proof%
Assume that $x^{4}+5x^{2}+5=5y^{2}$ for some positive integers $x$ and $y.$
Then, since $(2y+2)^{2}+(4y-1)^{2}=20y^{2}+5,$ it follows that 
\begin{equation*}
(2y+2)^{2}+(4y-1)^{2}=(2x^{2}+5)^{2}.
\end{equation*}%
Clearly, $d=(2y+2,4y-1)=1$ or $5.$ Assume that $d=1$. By the Pythagorean
theorem, there exist positive integers $a$ and $b$ with $(a,b)=1,$ $a$ and $%
b $ are opposite parity, such that 
\begin{equation*}
2x^{2}+5=a^{2}+b^{2},\text{ }2y+2=2ab,\text{ }4y-1=a^{2}-b^{2}.
\end{equation*}%
The latter two equations imply that 
\begin{equation}
-5=a^{2}-4ab-b^{2}.  \label{3.6}
\end{equation}%
Thus by Theorem \ref{t1}, we get $a=L_{3z+3}/2,$ $b=L_{3z}/2$ with
nonnegative even integer $z.$ On the other hand, from the equations $%
-5=a^{2}-4ab-b^{2}$ and $2x^{2}+5=a^{2}+b^{2},$ we readily obtain $%
x^{2}=a(a-2b).$ Since $(a,b)=1,$ it follows that, $r=(a,a-2b)=1$ or $2.$ If $%
r=1,$ then there exist coprime positive integers $u$ and $v$ such that $%
a=u^{2},$ $a-2b=v^{2}.$ Thus $L_{3z+3}=2a=2u^{2}$ and therefore $3z+3=6$ by (%
\ref{2.10}), which is impossible since $z$ is even. If $r=2,$ then $%
a=2u^{2}, $ $a-2b=2v^{2}.$ Thus $L_{3z+3}=4u^{2}=(2u)^{2}$ and therefore $%
3z+3=1$ or $3 $ by (\ref{2.9}). The first of these is impossible. And the
second implies that $z=0.$ Thus $a=2,$ $b=1.$ Since $2x^{2}+5=a^{2}+b^{2},$
it follows that $x=0,$ which is impossible since $x$ is positive. Assume
that $d=5.$ Then there exist positive integers $a$ and $b$ with $(a,b)=1,$ $%
a $ and $b$ are opposite parity, such that 
\begin{equation*}
2x^{2}+5=5a^{2}+5b^{2},\text{ }2y+2=10ab,\text{ }4y-1=5a^{2}-5b^{2}.
\end{equation*}%
The above first equation implies that $5|x$ and therefore $x=5t$ for some
positive integer $t.$ And the latter two equations imply that $%
-5=5a^{2}-20ab-5b^{2},$ i.e., $-1=a^{2}-4ab-b^{2}.$ and completing the
square gives $(a-2b)^{2}-5b^{2}=-1.$ Thus by Theorem \ref{t1}, we get $%
a=F_{3z+3}/2,$ $b=F_{3z}/2$ with positive odd integer $z.$ On the other
hand, by using $x=5t,$ from the equations $-5=5a^{2}-20ab-5b^{2}$ and $%
2x^{2}+5=5a^{2}+5b^{2},$ we obtain $5t^{2}=a(a-2b).$ Since $(a,b)=1,$
clearly, $(a,a-2b)=1$ or $2.$ Assume that $\ (a,a-2b)=1.$ This implies that
either $a=5u^{2},a-2b=v^{2}$ or $a=u^{2},a-2b=5v^{2}.$ If the first of these
is satisfied, then it is seen that $F_{3z+3}=10u^{2}$ and therefore $3z+3=0$
by (\ref{2.8}), which is impossible in positive integers. If the second is
satisfied, then it is seen that $F_{3z+3}=2u^{2}$ and therefore $3z+3=0,3$
or $6$ by (\ref{2.6}). But it is obvious that the cases $3z+3=0$ and $3z+3=3$
are impossible in positive integers. If $3z+3=6,$ then $z=1$ and therefore $%
a=2,b=1.$ Since $2x^{2}+5=5a^{2}+5b^{2},$ it follows that $x^{2}=10,$ which
is impossible. Assume that $\ (a,a-2b)=2.$ Then either $%
a=10u^{2},a-2b=2v^{2} $ or $a=2u^{2},a-2b=10v^{2}.$ If the first of these is
satisfied, then $F_{3z+3}=20u^{2}=5(2u)^{2}$ and therefore $3z+3=0$ or $5$
by (\ref{2.7}), which are impossible in positive integers. If the second is
satisfied, then $F_{3z+3}=4u^{2}=(2u)^{2}$ and therefore $3z+3=0,1,2$ or $12$
by (\ref{2.5}). But there does not any positive integer $z$ such that $%
3z+3=0,1$ or $2.$ If $3z+3=12,$ then we get $z=3$ and therefore $a=72,b=17.$
Since $2x^{2}+5=5a^{2}+5b^{2},$ it follows that $x^{2}=13680,$ which is
impossible. This completes the proof of Theorem \ref{t3.2}.%
\endproof%

We now state the following lemma without proof since its proof can be given
by induction method.

\begin{lemma}
\label{l2.3}If $n$ is even, then $V_{n}\equiv 2(\func{mod}$ $P^{2})$ and if $%
n$ is odd, then $V_{n}\equiv nP(\func{mod}$ $P^{2}).$
\end{lemma}

From Lemma \ref{l2.3} and identity (\ref{2.13}), we can give the following
corollary.

\begin{corollary}
\label{c2.2}$5|V_{n}$ if and only if $5|P$ and $n$ is odd.
\end{corollary}

The proof of the following lemma can be seen from identity (\ref{2.22}).

\begin{lemma}
\label{l2.5}If $P$ is odd and $r\geq 1,$ then $\left( \dfrac{P^{2}+3}{%
V_{2^{r}}}\right) =1.$
\end{lemma}

\begin{theorem}
\label{t3.3}If $P$ is odd, then the equation $V_{n}=5x^{2}$ has solutions
only if $n=1.$
\end{theorem}

\proof%
Assume that $V_{n}=5x^{2}.$ Then by Corollary \ref{c2.2}, it follows that $%
5|P$ and $n$ is odd. Assume that $n>3.$ Then we can write $n=4q+1$ or $%
n=4q+3 $ for some $q\geq 1.$ From this point on, we divide the proof into
two cases.

Case $1:$ Assume that $n=4q+1.$ Then we can write $n=4q+1=2(2^{k}a)+1$ for
some odd integer $a$ with $k\geq 1.$ And so by (\ref{2.4}), we get 
\begin{equation*}
V_{n}=V_{2.2^{k}a+1}\equiv -V_{1}(\func{mod}V_{2^{k}}),
\end{equation*}%
which implies that 
\begin{equation*}
5x^{2}\equiv -P(\func{mod}V_{2^{k}}).
\end{equation*}%
Therefore the Jacobi symbol $J=\left( \dfrac{-5P}{V_{2^{k}}}\right) =1.$
Assume that $P\equiv 5,7(\func{mod}$ $8).$ Since $V_{2^{k}}\equiv 2(\func{mod%
}$ $P)$ by Lemma \ref{l2.3}, it is seen that $V_{2^{k}}\equiv 2(\func{mod}$ $%
5).$ This shows that 
\begin{equation*}
\left( \dfrac{5}{V_{2^{k}}}\right) =\left( \frac{V_{2^{k}}}{5}\right)
=\left( \frac{2}{5}\right) =(-1)^{\frac{5^{2}-1}{8}}=-1
\end{equation*}%
and 
\begin{equation*}
\left( \dfrac{P}{V_{2^{k}}}\right) =(-1)^{\left( \frac{P-1}{2}\right) \left( 
\frac{V_{2^{r}}-1}{2}\right) }\left( \dfrac{V_{2^{k}}}{P}\right)
=(-1)^{\left( \frac{P-1}{2}\right) }\left( \frac{2}{P}\right) =(-1)^{\left( 
\frac{P-1}{2}\right) }(-1)^{\left( \frac{P^{2}-1}{8}\right) }=-1
\end{equation*}%
since $P\equiv 5,7(\func{mod}$ $8).$ Also we have $\left( \dfrac{-1}{%
V_{2^{k}}}\right) =-1$ by Lemma \ref{l2.2}. Hence we get $J=\left( \dfrac{-5P%
}{V_{2^{k}}}\right) =-1,$ which contradicts with the fact that $J=1.$ Assume
that $P\equiv 1,3(\func{mod}$ $8).$ If we write $n=4q+1=4(q+1)-3=2(2^{k}a)-3$
for some odd integer $a$ with $k\geq 1,$ then we get 
\begin{equation*}
V_{n}=V_{2.2^{k}a-3}\equiv -V_{-3}\equiv V_{3}(\func{mod}V_{2^{k}}),
\end{equation*}%
which implies that 
\begin{equation*}
5x^{2}\equiv V_{3}(\func{mod}V_{2^{k}})
\end{equation*}%
by (\ref{2.4}). This shows that $\left( \dfrac{5V_{3}}{V_{2^{k}}}\right) =1.$
Since $V_{2^{k}}\equiv 2(\func{mod}$ $P),$ we get $V_{2^{k}}\equiv 2(\func{%
mod}$ $5)$ by Lemma \ref{l2.3}. Moreover, $\left( \dfrac{P^{2}+3}{V_{2^{k}}}%
\right) =1$ by Lemma \ref{l2.5} and $V_{2^{k}}\equiv 3,7(\func{mod}$ $8)$ by
Lemma \ref{l2.2}. Then it follows that 
\begin{eqnarray*}
1 &=&\left( \frac{5V_{3}}{V_{2^{k}}}\right) =\left( \frac{5}{V_{2^{k}}}%
\right) \left( \frac{P}{V_{2^{k}}}\right) \left( \frac{P^{2}+3}{V_{2^{k}}}%
\right) =\left( \frac{V_{2^{k}}}{5}\right) (-1)^{^{\left( \frac{P-1}{2}%
\right) \left( \frac{V_{2^{r}-1}}{2}\right) }}\left( \frac{V_{2^{k}}}{P}%
\right) \\
&=&\left( \frac{2}{5}\right) (-1)^{\left( \frac{P-1}{2}\right) }\left( \frac{%
2}{P}\right) =(-1)(-1)^{\left( \frac{P-1}{2}\right) }(-1)^{\left( \frac{%
P^{2}-1}{8}\right) }=-1,
\end{eqnarray*}%
a contradiction.

Case $2:$ Assume that $n=4q+3.$ We can write $n=4q+3=2(2^{k}a)+3$ for some
odd integer $a$ with $k\geq 1.$ And so by (\ref{2.4}), we get 
\begin{equation*}
V_{n}=V_{2.2^{k}a+3}\equiv -V_{3}(\func{mod}\text{ }V_{2^{k}})
\end{equation*}%
i.e.,%
\begin{equation*}
5x^{2}=-V_{3}(\func{mod}\text{ }V_{2^{k}}).
\end{equation*}%
This shows that $J=\left( \dfrac{-5V_{3}}{V_{2^{k}}}\right) =1.$ Assume that 
$P\equiv 5,7(\func{mod}$ $8).$ Since $V_{2^{k}}\equiv 2(\func{mod}$ $P)$ by
Lemma \ref{l2.3}, it is seen that $V_{2^{k}}\equiv 2(\func{mod}$ $5).$ Also
we have $\left( \dfrac{-1}{V_{2^{k}}}\right) =-1$ by Lemma \ref{l2.2} and $%
\left( \dfrac{P^{2}+3}{V_{2^{k}}}\right) =1$ by Lemma \ref{l2.5}. Hence we
get 
\begin{eqnarray*}
\left( \frac{-5V_{3}}{V_{2^{k}}}\right) &=&\left( \dfrac{-1}{V_{2^{k}}}%
\right) \left( \frac{5}{V_{2^{k}}}\right) \left( \frac{V_{3}}{V_{2^{k}}}%
\right) =\left( \dfrac{-1}{V_{2^{k}}}\right) \left( \frac{5}{V_{2^{k}}}%
\right) \left( \frac{P}{V_{2^{k}}}\right) \left( \frac{P^{2}+3}{V_{2^{k}}}%
\right) \\
&=&(-1)\left( \frac{V_{2^{k}}}{5}\right) (-1)^{\left( \frac{P-1}{2}\right)
\left( \frac{V_{2^{r}-1}}{2}\right) }\left( \frac{V_{2^{k}}}{P}\right)
=(-1)\left( \frac{2}{5}\right) (-1)^{\left( \frac{P-1}{2}\right) }\left( 
\frac{2}{P}\right) \\
&=&(-1)(-1)(-1)^{\left( \frac{P-1}{2}\right) }(-1)^{\left( \frac{P^{2}-1}{8}%
\right) }=-1
\end{eqnarray*}%
since $P\equiv 5,7(\func{mod}$ $8).$ This contradicts with the fact that $%
J=1.$ Assume that $P\equiv 1,3(\func{mod}$ $8).$ If we write $%
n=4q+3=4(q+1)-1=2(2^{k}a)-1$ for some odd integer $a$ with $k\geq 1,$ then
we get 
\begin{equation*}
V_{n}=V_{2.2^{k}a-1}\equiv -V_{-1}\equiv V_{1}(\func{mod}\text{ }V_{2^{k}})
\end{equation*}%
i.e.,%
\begin{equation*}
5x^{2}=P(\func{mod}\text{ }V_{2^{k}}).
\end{equation*}%
This shows that $\left( \dfrac{5P}{V_{2^{k}}}\right) =1.$ Since $%
V_{2^{k}}\equiv 2(\func{mod}$ $P),$ we get $V_{2^{k}}\equiv 2(\func{mod}$ $%
5) $ by Lemma \ref{l2.3}. Then it follows that 
\begin{eqnarray*}
1 &=&\left( \frac{5P}{V_{2^{k}}}\right) =\left( \frac{5}{V_{2^{k}}}\right)
\left( \frac{P}{V_{2^{k}}}\right) =\left( \frac{V_{2^{k}}}{5}\right)
(-1)^{\left( \frac{P-1}{2}\right) \left( \frac{V_{2^{r}-1}}{2}\right)
}\left( \frac{V_{2^{k}}}{P}\right) \\
&=&\left( \frac{2}{5}\right) (-1)^{\left( \frac{P-1}{2}\right) }\left( \frac{%
2}{P}\right) =(-1)(-1)^{\left( \frac{P-1}{2}\right) }(-1)^{\left( \frac{%
P^{2}-1}{8}\right) }=-1,
\end{eqnarray*}%
a contradiction. We conclude that $n=1$ or $n=3.$ If $n=3,$ then $%
V_{3}=P(P^{2}+3)=5x^{2}.$ Since $5|P,$ it follows that $%
(P/5)(P^{2}+3)=x^{2}. $ Clearly, $d=(P/5,P^{2}+3)=1$ or $3.$ Assume that $%
d=1.$ This implies that $P=5a^{2}$ and $P^{2}+3=b^{2}$ for some positive
integers $a$ and $b.$ Since $5|P,$ we get $b^{2}\equiv 3(\func{mod}$ $5),$
which is impossible. Assume that $d=3.$ Then we get $P=15a^{2}$ and $%
P^{2}+3=3b^{2}$ for some positive integers $a$ and $b.$ It is seen from $%
P^{2}+3=3b^{2}$ that $3|P$ and therefore $P=3c$ for some positive integer $%
c. $ Hence we obtain the Pell equation $b^{2}-3c^{2}=1.$ It is well known
that all positive integer solutions of this equation are given by $%
(b,c)=(v_{m}(4,-1)/2,u_{m}(4,-1))$ with $m\geq 1.$ On the other hand, if we
substitute the value $P=15a^{2}$ into $P=3c,$ we get $c=5a^{2}.$ So we are
interested in whether the equation $5\square =u_{m}(4-1)$ has a solution.
Assume that the equation $5\square =u_{m}(4-1)$ has a solution. Since $%
5|u_{3},$ it can be seen that if $5|u_{m},$ then $3|m$ and therefore $m=3r$
for some positive integer $r.$ Thus from (\ref{10})\ we get $%
u_{m}=u_{3r}=u_{r}\left( (P^{2}-4)u_{r}^{2}+3\right) =u_{r}(12u_{r}^{2}+3).$
Clearly, $(u_{r},12u_{r}^{2}+3)=1$ or $3.$ Assume that $%
(u_{r},12u_{r}^{2}+3)=1.$ This implies that either $u_{r}=a^{2},$ $%
12u_{r}^{2}+3=5b^{2}$ or $u_{r}=5a^{2},$ $12u_{r}^{2}+3=b^{2}\ $for some
positive integers $a$ and $b.$ But both of the previous equations are
impossible since $b^{2}\equiv 3(\func{mod}$ $5).$ Assume that $%
(u_{r},12u_{r}^{2}+3)=3.$ Then either 
\begin{equation}
u_{r}=3a^{2},\text{ }12u_{r}^{2}+3=15b^{2}  \label{11}
\end{equation}%
or%
\begin{equation}
u_{r}=15a^{2},\text{ }12u_{r}^{2}+3=3b^{2}.  \label{12}
\end{equation}%
Assume that (\ref{11}) is satisfied. A simple computation shows that $%
(2u_{r})^{2}-5b^{2}=-1.$ Thus by Theorem \ref{t2.8}, we obtain $%
2u_{r}=L_{3z}/2$ for some positive odd integer $z.$ Substituting the value $%
u_{r}=3a^{2}$ into the previous equation gives $3u_{r}=L_{3z}/4,$ i.e., $%
L_{2}u_{r}=L_{3z}/4.$ This implies that $L_{2}|L_{3z}.$ Then by (\ref{2.17}%
), we get $2|3z,$ which is impossible since $z$ is odd. Assume that (\ref{12}%
) is satisfied. It is easily seen that $(2u_{r})^{2}+1=b^{2},$ that is, $%
b^{2}-(2u_{r})^{2}=1,$ implying that $u_{r}=0.$ This is impossible since $r$
is a positive integer. So $n=3$ can not be a solution. If $\ n=1,$ then $%
V_{1}=P=5x^{2}.$ It is obvious that this is a solution. This completes the
proof of Theorem \ref{t3.3}.%
\endproof%

\begin{theorem}
\label{t3.4}There is no integer $x$ such that $V_{n}=5V_{m}x^{2}.$
\end{theorem}

\proof%
Assume that $V_{n}=5V_{m}x^{2}.$ Then by Corollary \ref{c2.2}, it follows
that $5|P$ and $n$ is odd. Moreover, since $V_{m}|V_{n},$ there exists an
odd integer $t$ such that $n=mt$ by (\ref{2.17}). Thus $m$ is odd. Therefore
we have $V_{n}\equiv nP(\func{mod}$ $P^{2})$ and $V_{m}\equiv mP(\func{mod}$ 
$P^{2})$ by Lemma \ref{l2.3}. This shows that $nP\equiv 5mPx^{2}(\func{mod}$ 
$P^{2}),$ i.e., $n\equiv 5mx^{2}(\func{mod}$ $P).$ Since $5|P,$ it follows
that $5|n.$ Also since $n=mt,$ first, assume that $5|t.$ Then $t=5s$ for
some positive odd integer $s$ and therefore $n=mt=5ms.$ By (\ref{2.16}), we
readily obtain $V_{n}=V_{5ms}=V_{ms}(V_{ms}^{4}+5V_{ms}^{2}+5).$ Since $ms$
is odd and $5|P,$ it follows that $5|V_{ms}$ by Corollary \ref{c2.2} and
therefore $(V_{ms}/V_{m})((V_{ms}^{4}+5V_{ms}^{2}+5)/5)=x^{2}.$ Clearly, $%
(V_{ms}/V_{m},(V_{ms}^{4}+5V_{ms}^{2}+5)/5)=1.$ This implies that $%
V_{ms}=V_{m}a^{2}$ and $V_{ms}^{4}+5V_{ms}^{2}+5=5b^{2}$ for some positive
integers $a$ and $b.$ Then by Theorem \ref{t3.2}, we get $V_{ms}=0,$ which
is a contradiction. Now assume that $5\nmid t.$ Since $n=mt$ and $5|n,$ it
is seen that $5|m.$ Then we can write $m=5^{r}a$ with $5\nmid a$ and $r\geq
1.$ By (\ref{2.101}), we obtain $V_{m}=V_{5^{r}a}=5V_{5^{r-1}a}(5a_{1}+1)$
for some positive integer $a_{1}.$ And thus we conclude that $%
V_{m}=V_{5^{r}a}=5^{r}V_{a}(5a_{1}+1)(5a_{2}+1)...(5a_{r}+1)$ for some
positive integers $a_{i}$ with $1\leq i\leq r.$ Let $%
A=(5a_{1}+1)(5a_{2}+1)...(5a_{r}+1).$ It is obvious that $5\nmid A.$ Thus we
have $V_{m}=5^{r}V_{a}A.$ In a similar manner, we see that $%
V_{n}=V_{5^{r}at}=5^{r}V_{at}(5b_{1}+1)(5b_{2}+1)...(5b_{r}+1)$ for some
positive integers $b_{j}$ with $1\leq j\leq r.$ Let $%
B=(5b_{1}+1)(5b_{2}+1)...(5b_{r}+1).$ It is obvious that $5\nmid B.$ Thus we
have $V_{n}=5^{r}V_{at}B.$ This shows that $5^{r}V_{at}B=5.5^{r}V_{a}Ax^{2},$
i.e., $V_{at}B=5V_{a}Ax^{2}.$ By Lemma \ref{l2.3} and Corollary \ref{c2.2},
it is seen that $atPB\equiv 5aPAx^{2}(\func{mod}$ $P^{2})$ and therefore we
get $atB\equiv 5aAx^{2}(\func{mod}$ $P).$ Since $5|P,$ it follows that $%
5|atB.$ But this is impossible since $5\nmid a,5\nmid t,$ and $5\nmid B.$
This completes the proof of Theorem \ref{t3.4}.%
\endproof%

The following lemma can be proved by using Theorem \ref{t2.1}.

\begin{lemma}
\label{l2.4}%
\begin{equation*}
5|U_{n}\Leftrightarrow \left\{ 
\begin{array}{c}
2|n\text{ if }5|P, \\ 
3|n\text{ if }P^{2}\equiv -1(\func{mod}\text{ }5), \\ 
5|n\text{ if }P^{2}\equiv 1(\func{mod}\text{ }5),%
\end{array}%
\right.
\end{equation*}%
and 
\begin{equation*}
3|U_{n}\Leftrightarrow \left\{ 
\begin{array}{c}
2|n\text{ if }3|P, \\ 
4|n\text{ if }3\nmid P.%
\end{array}%
\right.
\end{equation*}
\end{lemma}

\begin{theorem}
\label{t3.6}If $P$ is odd and $5|P,$ then the equation $U_{n}=5x^{2}$ has
the solution $n=2,P=5\square .$ If \ $P^{2}\equiv 1(\func{mod}$ $5),$ then
the equation $U_{n}=5x^{2}$ has the solution $n=5,$ $P=1.$ If $P$ is odd and 
$P^{2}\equiv -1(\func{mod}$ $5),$ then the equation $U_{n}=5x^{2}$ has no
solutions.
\end{theorem}

\proof%
Assume that $5|P$ and $P$ is odd. Since $5|U_{n},$ it follows that $n$ is
even by Lemma \ref{l2.4}. Then $n=2t$ for some positive integer $t.$ By (\ref%
{2.11}), we get $U_{n}=U_{2t}=U_{t}V_{t}=5x^{2}.$ Clearly, $(U_{t},V_{t})=1$
or $2$ by (\ref{2.18}). Let $(U_{t},V_{t})=1.$ This implies that either 
\begin{equation}
U_{t}=a^{2},\text{ }V_{t}=5b^{2}  \label{3.7}
\end{equation}%
or%
\begin{equation}
U_{t}=5a^{2},\text{ }V_{t}=b^{2}  \label{3.8}
\end{equation}%
for some positive integers $a$ and $b.$ Assume that (\ref{3.7}) is
satisfied. Since $5|V_{t},$ it follows that $t$ is an odd integer by
Corollary \ref{c2.2}. Assume that $t>1.$ Then $t=4q\pm 1$ for some $q>1.$ We
can write $t=4q\pm 1=2.2^{k}u\pm 1$ for some odd integer $u$ with $k\geq 1.$
And so by (\ref{2.3}), we get 
\begin{equation*}
U_{t}=U_{2.2^{k}u\pm 1}\equiv -U_{\pm 1}(\func{mod}\text{ }V_{2^{k}})
\end{equation*}%
which implies that 
\begin{equation*}
a^{2}\equiv -1(\func{mod}\text{ }V_{2^{k}}).
\end{equation*}%
This shows that $1=\left( \frac{-1}{V_{2^{k}}}\right) .$ But this is
impossible since $\left( \frac{-1}{V_{2^{k}}}\right) =-1$ by Lemma \ref{l2.2}%
. Thus $t=1$ and therefore $n=2.$ Then $P=5\square $ is a solution. Assume
that (\ref{3.8}) is satisfied. Since $5|U_{t},$ it follows that $t$ is even
by Lemma \ref{l2.4}. Thus $t=2r$ for some positive integer $r.$ By using (%
\ref{2.12}), we get $V_{2r}=V_{r}^{2}\pm 2=b^{2},$ which is impossible. Thus 
$t=1$ and therefore $n=2.$ Let $d=2.$ This implies that either 
\begin{equation}
U_{t}=10a^{2},\text{ }V_{t}=2b^{2}  \label{3.9}
\end{equation}%
or%
\begin{equation}
U_{t}=2a^{2},\text{ }V_{t}=10b^{2}  \label{3.10}
\end{equation}%
for some positive integers $a$ and $b.$ Assume that (\ref{3.9}) is
satisfied. By Theorem \ref{t2.5}, we have $t=6$ and $P=5.$ But this is
impossible since there does not exist any integer $a$ such that $%
U_{6}=3640=10a^{2}.$ Assume that (\ref{3.10}) is satisfied. Since $5|V_{t},$
it follows that $t$ is an odd integer by Corollary \ref{c2.2}. If $t=1,$
then $U_{1}=1=2a^{2},$ which is impossible. Assume that $t>1.$ Then $t=4q\pm
1$ for some $q>1.$ And so by (\ref{2.1}), we get%
\begin{equation*}
U_{t}=U_{2.2q\pm 1}\equiv U_{\pm 1}(\func{mod}\text{ }U_{2}),
\end{equation*}%
implying that 
\begin{equation*}
2a^{2}\equiv 1(\func{mod}\text{ }P).
\end{equation*}%
Since $5|P,$ the above congruence becomes%
\begin{equation*}
2a^{2}\equiv 1(\func{mod}\text{ }5),
\end{equation*}%
which is impossible since $\left( \dfrac{2}{5}\right) =-1.$ The proof is
completed for the case when $5|P$ and $P$ is odd.

Assume that $P^{2}\equiv 1(\func{mod}$ $5).$ Since $5|U_{n},$ it follows
that $5|n$ by Lemma \ref{l2.4}. Thus $n=5t$ for some positive integer $t.$
Since $P^{2}\equiv 1(\func{mod}$ $5),$ it is obvious that $5|P^{2}+4$ and
therefore there exists a positive integer $A$ such that $P^{2}+4=5A.$ By (%
\ref{2.15}), we get $U_{n}=U_{5t}=U_{t}\left( (P^{2}+4)^{2}U_{t}^{4}\pm
5(P^{2}+4)U_{t}^{2}+5\right) .$ Substituting $P^{2}+4=5A$ into the previous
equation gives $U_{n}=U_{5t}=5U_{t}(5A^{2}U_{t}^{4}\pm 5AU_{t}^{2}+1).$ Let $%
B=A^{2}U_{t}^{4}\pm AU_{t}^{2}.$ Then we get%
\begin{equation*}
U_{n}=U_{5t}=5U_{t}(5B+1)=5x^{2}
\end{equation*}%
i.e.,%
\begin{equation*}
U_{t}(5B+1)=x^{2}.
\end{equation*}%
It can be seen that $(U_{t},5B+1)=1.$ This shows that $U_{t}=a^{2}$ and $%
5B+1=b^{2}$ for some positive integers $a$ and $b.$ By Theorem \ref{t2.3},
we get $t\leq 2$ or $t=12$ and $P=1.$ If $t=1,$ then $n=5$ and therefore we
get $U_{5}=P^{4}+3P^{2}+1=5x^{2}.$ By Theorem \ref{t3.1}, it follows that $%
P=1.$ So the equation $U_{n}=5x^{2}$ has the solution $n=5$ and $P=1.$ If $%
t=2,$ then $n=10$ and therefore we obtain $U_{10}=5x^{2},$ implying that $%
U_{5}V_{5}=5x^{2}$ by (\ref{2.11}). Since $5|U_{5},$ it follows that $%
(U_{5}/5)V_{5}=x^{2}.$ By (\ref{2.18}), clearly, $(U_{5}/5,V_{5})=1.$ This
implies that $U_{5}=5a^{2},$ $V_{5}=b^{2},$ which is impossible by Theorem %
\ref{t2.4}. If $t=12$ and $P=1,$ then it follows that $n=60.$ Thus we obtain 
$U_{60}=5x^{2},$ which is impossible by (\ref{2.7}). The proof is completed
for the case when $P^{2}\equiv 1(\func{mod}$ $5).$

Assume that $P^{2}\equiv -1(\func{mod}$ $5)$ and $P$ is odd. Since $5|U_{n},$
it follows that $3|n$ by Lemma \ref{l2.4} and therefore $n=3m$ for some
positive integer $m.$ Assume that $m$ is even. Then $m=2s$ for some positive
integer $s$ and therefore $n=6s.$ Thus by (\ref{2.11}), we get $%
U_{n}=U_{6s}=U_{3s}V_{3s}=5x^{2}.$ By (\ref{2.18}), clearly, $%
(U_{3s},V_{3s})=2.$ Then either 
\begin{equation}
U_{3s}=10a^{2},\text{ }V_{3s}=2b^{2}  \label{3.14}
\end{equation}%
or%
\begin{equation}
U_{3s}=2a^{2},\text{ }V_{3s}=10b^{2}  \label{3.15}
\end{equation}%
for some positive integer $a$ and $b.$ Assume that (\ref{3.14}) is
satisfied. By Theorem \ref{t2.5}, it follows that $3s=6$ and $P=1,5.$ But
this is impossible since $P^{2}\equiv -1(\func{mod}$ $5).$ Assume that (\ref%
{3.15}) is satisfied. Since $5|V_{3s},$ it follows that $5|P$ by Corollary %
\ref{c2.2}. But this contradicts with the fact that $P^{2}\equiv -1(\func{mod%
}$ $5).$ Now assume that $m$ is odd. Then by (\ref{2.14}), we get $%
U_{n}=U_{3m}=U_{m}\left( (P^{2}+4)U_{m}^{2}-3\right) .$ Clearly, $%
(U_{m},(P^{2}+4)U_{m}^{2}-3)=1$ or $3.$ Since $m$ is odd, it follows that $%
3\nmid U_{m}$ by Lemma \ref{l2.4} and therefore $%
(U_{m},(P^{2}+4)U_{m}^{2}-3)=1.$ Then 
\begin{equation}
U_{m}=5a^{2},\text{ }(P^{2}+4)U_{m}^{2}-3=b^{2}  \label{3.16}
\end{equation}%
or%
\begin{equation}
U_{m}=a^{2},\text{ }(P^{2}+4)U_{m}^{2}-3=5b^{2}  \label{3.17}
\end{equation}%
for some positive integers $a$ and $b.$ Assume that (\ref{3.16}) is
satisfied. Since $m$ is odd, we obtain $V_{m}^{2}+1=b^{2}$ by (\ref{2.13}).
This shows that $V_{m}=0,$ which is impossible. Assume that (\ref{3.17}) is
satisfied. Since $m$ and $P$ is odd, it follows that $m=1$ by Theorem \ref%
{t2.3}. If $m=1,$ then $n=3$ and therefore $P^{2}+1=5y^{2},$ which is
impossible since we get $y^{2}\equiv 2(\func{mod}$ $8)$ in this case. This
completes the proof of Theorem \ref{t3.6}.%
\endproof%

Since the proof of the following lemma can be given by induction method, we
omit its proof.

\begin{lemma}
\label{l2.7}If $n$ is even, then $U_{n}\equiv \dfrac{n}{2}P(\func{mod}$ $%
P^{2})$ and if $n$ is odd, then $U_{n}\equiv 1(\func{mod}$ $P^{2}).$
\end{lemma}

\begin{theorem}
\label{t3.7}The equation $U_{n}=5U_{m}x^{2}$ has no solutions when $%
P^{2}\equiv 1(\func{mod}$ $5).$ If $P$ is odd or $4|P,$ then the equation $%
U_{n}=5U_{m}x^{2}$ has no solutions when $P^{2}\equiv -1(\func{mod}$ $5)$
and $n$ is odd. If $n$ is even and $P$ is odd, then the equation $%
U_{n}=5U_{m}x^{2}$ has no solutions when $P^{2}\equiv -1(\func{mod}$ $5).$
If $P$ is odd and $5|P,$ then the equation $U_{n}=5U_{m}x^{2}$ has no
solutions.
\end{theorem}

\proof%
Assume that $U_{n}=5U_{m}x^{2}$ for some positive integer $x.$ If $m=1,$
then $U_{n}=5x^{2}$ which has solutions only if $n=2$ by Theorem \ref{t3.6}%
.\ So assume that $m>1.$ Since $U_{m}|U_{n},$ it follows that $m|n$ by (\ref%
{2.01}). Thus $n=mt$ for some positive integer $t.$ Since $n\neq m,$ we have 
$t>1.$

Assume that $P^{2}\equiv 1(\func{mod}$ $5).$ It is obvious that $5|P^{2}+4.$
Since $5|U_{n},$ it follows that $5|n$ by Lemma \ref{l2.4}. Now we divide
the proof into two cases.

Case $1:$ Assume that $5|t.$ Then $t=5s$ for some positive integer $s$ and
therefore $n=mt=5ms.$ By (\ref{2.15}), we obtain 
\begin{equation}
U_{n}=U_{5ms}=U_{ms}\left( (P^{2}+4)^{2}U_{ms}^{4}\pm
5(P^{2}+4)U_{ms}^{2}+5\right) =5U_{m}x^{2}.  \label{3.20}
\end{equation}%
It is easily seen that $5|(P^{2}+4)^{2}U_{ms}^{4}\pm 5(P^{2}+4)U_{ms}^{2}+5.$
Also we have $(P^{2}+4)^{2}U_{ms}^{4}\pm
5(P^{2}+4)U_{ms}^{2}+5=V_{ms}^{4}\pm 3V_{ms}^{2}+1$ by (\ref{2.13}). So
rearranging the equation (\ref{3.20}) gives 
\begin{equation*}
x^{2}=(U_{ms}/U_{m})\left( (V_{ms}^{4}\pm 3V_{ms}^{2}+1)/5\right) .
\end{equation*}%
Clearly, $(U_{ms}/U_{m},(V_{ms}^{4}\pm 3V_{ms}^{2}+1)/5)=1.$ This implies
that $U_{ms}=U_{m}a^{2}$ and $V_{ms}^{4}\pm 3V_{ms}^{2}+1=5b^{2}$ for some
positive integers $a$ and $b.$ Thus by Theorem \ref{t3.1}, we get $V_{ms}=1$
or $V_{ms}=2.$ The first of these is impossible. If the second is satisfied,
then $ms=0,$ which is a contradiction since $m>1.$

Case $2:$ Assume that $5\nmid t.$ Since $5|n,$ it follows that $5|m.$ Then
we can write $m=5^{r}a$ with $5\nmid a$ and $r\geq 1.$ Since $5|P^{2}+4,$ it
can be seen by (\ref{2.15}) that $U_{m}=U_{5^{r}a}=5U_{5^{r-1}a}(5a_{1}+1)$
for some positive integer $a_{1}.$ And thus we conclude that $%
U_{m}=U_{5^{r}a}=5^{r}U_{a}(5a_{1}+1)(5a_{2}+1)...(5a_{r}+1)$ for some
positive integers $a_{i}$ with $1\leq i\leq r.$ Let $%
A=(5a_{1}+1)(5a_{2}+1)...(5a_{r}+1).$ It is obvious that $5\nmid A$ and we
have $U_{m}=5^{r}U_{a}A.$ In a similar manner, we get $%
U_{n}=U_{5^{r}at}=5^{r}U_{at}(5b_{1}+1)(5b_{2}+1)...(5b_{r}+1)$ for some
positive integers $b_{j}$ with $1\leq j\leq r.$ Let $%
B=(5b_{1}+1)(5b_{2}+1)...(5b_{r}+1).$ It is obvious that $5\nmid B.$ Thus we
have $U_{n}=5^{r}U_{at}B.$ Substituting the new values of $U_{n}$ and $U_{m}$
into $U_{n}=5U_{m}x^{2}$ gives 
\begin{equation*}
5^{r}U_{at}B=5.5^{r}U_{a}Ax^{2}.
\end{equation*}%
This shows that 
\begin{equation*}
U_{at}B=5U_{a}Ax^{2}.
\end{equation*}%
Since $5\nmid B,$ it follows that $5|U_{at},$ implying that $5|at$ by Lemma %
\ref{l2.4}. This contradicts with the fact that $5\nmid a$ and $5\nmid t.$

Assume that $P^{2}\equiv -1(\func{mod}$ $5)$ and $n$ is odd. Then, both $m$
and $t$ are odd. Thus we can write $t=4q\pm 1$ for some $q\geq 1.$ And so by
(\ref{2.1}), we get 
\begin{equation*}
U_{n}=U_{(4q\pm 1)m}=U_{2.2mq\pm m}\equiv U_{m}(\func{mod}\text{ }U_{2m}).
\end{equation*}%
This shows that 
\begin{equation*}
5U_{m}x^{2}\equiv U_{m}(\func{mod}\text{ }U_{2m}).
\end{equation*}%
By using (\ref{2.11}), we obtain 
\begin{equation*}
5x^{2}\equiv 1(\func{mod}\text{ }V_{m}).
\end{equation*}%
Since $m$ is odd, it follows that $P|V_{m}$ by Lemma \ref{l2.3}. Then the
above congruence becomes 
\begin{equation}
5x^{2}\equiv 1(\func{mod}\text{ }P).  \label{4}
\end{equation}%
Assume that $P$ is odd. Then (\ref{4}) implies that $J=\left( \dfrac{5}{P}%
\right) =1.$ Since $P^{2}\equiv -1(\func{mod}$ $5),$ it can be seen that $%
P\equiv \pm 2(\func{mod}$ $5).$ Hence we get 
\begin{equation*}
1=\left( \frac{5}{P}\right) =\left( \frac{P}{5}\right) =\left( \frac{\pm 2}{5%
}\right) =-1,
\end{equation*}%
a contradiction. Now assume that $P$ is even. If $8|P,$ then it follows from
(\ref{4}) that $5x^{2}\equiv 1(\func{mod}$ $8),$ which is impossible since
we get $x^{2}\equiv 5(\func{mod}$ $8)$ in this case. If $4|P$ and $8\nmid P,$
then from (\ref{4}), we get 
\begin{equation*}
5x^{2}\equiv 1(\func{mod}\text{ }P/4).
\end{equation*}%
This shows that $\left( \dfrac{5}{P/4}\right) =1.$ Since $P^{2}\equiv -1(%
\func{mod}$ $5),$ it can be seen that $P/4\equiv \pm 2(\func{mod}$ $5).$
Hence we get 
\begin{equation*}
1=\left( \frac{5}{P/4}\right) =\left( \frac{P/4}{5}\right) =\left( \frac{\pm
2}{5}\right) =-1,
\end{equation*}%
a contradiction.

Now assume that $P^{2}\equiv -1(\func{mod}$ $5),$ $P$ is odd, and $n$ is
even. Since $n=mt,$ we divide the proof into two cases.

Case $1:$ Assume that $t$ is even. Then $t=2s$ for some positive integer $s.$
Thus we get $5x^{2}=U_{n}/U_{m}=U_{2ms}/U_{m}=(U_{ms}/U_{m})V_{ms}.$
Clearly, $d=(U_{ms}/U_{m},V_{ms})=1$ or $2$ by (\ref{2.18}). Let $d=1.$ Then
either 
\begin{equation}
U_{ms}=U_{m}a^{2}\text{ and }V_{ms}=5b^{2}  \label{3.21}
\end{equation}%
or $\ $%
\begin{equation}
U_{ms}=5U_{m}a^{2}\text{ and }V_{ms}=b^{2}.  \label{3.22}
\end{equation}%
Assume that (\ref{3.21}) is satisfied. Since $5|V_{ms},$ it follows that $%
5|P $ by Corollary \ref{c2.2}. This contradicts with the fact that $%
P^{2}\equiv -1(\func{mod}$ $5).$ Assume that (\ref{3.22}) is satisfied. By
Theorem \ref{t2.4}, we get $ms=3$ and $P=3.$ Since $m>1,$ it follows that $%
m=3.$ This is impossible since we get $1=5a^{2}$ in this case.

Let $d=2.$ This implies that either 
\begin{equation}
U_{ms}=2U_{m}a^{2}\text{ and }V_{ms}=10b^{2}  \label{3.23}
\end{equation}%
or 
\begin{equation}
U_{ms}=10U_{m}a^{2}\text{ and }V_{ms}=2b^{2}.  \label{3.24}
\end{equation}%
Assume that (\ref{3.23}) is satisfied. Since $5|V_{ms},$ it follows that $%
5|P $ by Corollary \ref{c2.2}. This contradicts with the fact that $%
P^{2}\equiv -1(\func{mod}$ $5).$ Assume that (\ref{3.24}) is satisfied. By
Theorem \ref{t2.5}, we get $ms=6$ and $P=1,5.$ But this is impossible since $%
P^{2}\equiv -1(\func{mod}$ $5).$

Case $2:$ Assume that $t$ is odd. Since $n$ is even, it follows that $m$ is
even. Then there exists a positive integer $s$ such that $m=2s.$ Thus we
readily obtain $5x^{2}=(U_{st}/U_{s})(V_{st}/V_{s}).$ Clearly, $%
d=(U_{st}/U_{s},V_{st}/V_{s})=1$ or $2$ by (\ref{2.18}). Let $d=1.$ Then
either $U_{st}=U_{s}a^{2}$ and $V_{st}=5V_{s}b^{2}$ or $U_{st}=5U_{s}a^{2}$
and $V_{st}=V_{s}b^{2}$ for some positive integers $a$ and $b.$ The first of
these is impossible by Theorem \ref{t3.4}. If the second is satisfied, then
we get $st=s$ by Theorem \ref{T2.6}. But this impossible since there does
not exist any integer $a$ such that $1=5a^{2}.$ Let $d=2.$ This implies that
either $U_{st}=2U_{s}a^{2}$ and $V_{st}=10V_{s}b^{2}$ or $%
U_{st}=10U_{s}a^{2} $ and $V_{st}=2V_{s}b^{2}$ for some positive integers $a$
and $b.$ If the first of these is satisfied, then it follows that $5|V_{st}.$
This implies that $5|P$ by Corollary \ref{c2.2}, which contradicts with the
fact that $P^{2}\equiv -1(\func{mod}$ $5).$ The second is impossible by
Theorem \ref{T2.7}.

Now assume that $5|P$ and $P$ is odd. Since $5|U_{n},$ it follows that $n$
is even by Lemma \ref{l2.4}. Moreover, since $U_{m}|U_{n},$ there exists an
integer $t$ such that $n=mt$ by (\ref{2.01}). Assume that $t$ is even. Then $%
t=2s$ for some positive integer $s.$ By (\ref{2.11}), we get $%
U_{n}=U_{2ms}=U_{ms}V_{ms}=5U_{m}x^{2},$ implying that $%
(U_{ms}/U_{m})V_{ms}=5x^{2}.$ Clearly, $(U_{ms}/U_{m},V_{ms})=1$ or $2$ by (%
\ref{2.18}). If $(U_{ms}/U_{m},V_{ms})=1,$ then%
\begin{equation}
U_{ms}=U_{m}a^{2},\text{ }V_{ms}=5b^{2}  \label{1}
\end{equation}%
or%
\begin{equation}
U_{ms}=5U_{m}a^{2},\text{ }V_{ms}=b^{2}  \label{2}
\end{equation}%
for some positive integers $a$ and $b.$ Assume that (\ref{1}) is satisfied.
Then by Theorem \ref{t3.3}, we get $ms=1.$ This contradicts with the fact
that $m>1.$ Assume that (\ref{2}) is satisfied. Then by Theorem \ref{t2.4},
we have $ms=3$ and $P=1$ or $ms=3$ and $P=3.$ But both of these are
impossible since $5|P.$ If $(U_{ms}/U_{m},V_{ms})=2,$ then

\begin{equation}
U_{ms}=2U_{m}a^{2},\text{ }V_{ms}=10b^{2}  \label{3}
\end{equation}%
or%
\begin{equation}
U_{ms}=10U_{m}a^{2},\text{ }V_{ms}=2b^{2}  \label{5}
\end{equation}%
for some positive integers $a$ and $b.$ Assume that (\ref{3}) is satisfied.
Then by Teorem \ref{T2.9}, we get $ms=12,$ $m=6,$ $P=5.$ On the other hand,
since $5|V_{ms},$ it follows by Corollary \ref{c2.2} that $5|P$ and $ms$ is
odd. This is a contradiction since $ms=12.$ Assume that (\ref{5}) is
satisfied. Then by Theorem \ref{t2.5}, we have $ms=6$ and $P=5.$ Since $m>1,$
it is seen that $m=2,3$ or $6.$ If $m=2,$ then $%
U_{6}=3640=10U_{2}x^{2}=50x^{2},$ i.e., $364=5x^{2},$ which is impossible.
If $m=3,$ then $U_{6}=3640=10U_{3}x^{2}=260x^{2},$ i.e., $14=x^{2},$ which
is impossible. If $m=6,$ then there does not exists any integer $x$ such
that $1=5x^{2}.$

Now assume that $t$ is odd. Since $n=mt$ and $n$ is even, it follows that $m$
is even. Therefore we have $U_{n}\equiv (n/2)P(\func{mod}$ $P^{2})$ and $%
U_{m}\equiv (m/2)P(\func{mod}$ $P^{2})$ by Lemma \ref{l2.7}. This shows that 
$(n/2)P\equiv 5(m/2)Px^{2}(\func{mod}$ $P^{2}),$ i.e., $(n/2)\equiv
5(m/2)x^{2}(\func{mod}$ $P).$ Since $5|P,$ it is obvious that $5|n.$ Now we
divide the proof into two cases.

Case $1:$ Assume that $5|t.$ Then $t=5s$ for some positive integer $s$ and
therefore $n=mt=5ms.$ By (\ref{2.15}), we obtain 
\begin{equation}
U_{n}=U_{5ms}=U_{ms}\left(
(P^{2}+4)^{2}U_{ms}^{4}+5(P^{2}+4)U_{ms}^{2}+5\right) =5U_{m}x^{2}.
\label{6}
\end{equation}%
Since $ms$ is even and $5|P,$ it is seen that $5|U_{ms}$ by \ref{l2.4}. Also
we have $%
(P^{2}+4)^{2}U_{ms}^{4}+5(P^{2}+4)U_{ms}^{2}+5=V_{ms}^{4}-3V_{ms}^{2}+1$ by (%
\ref{2.13}). So rearranging the equation (\ref{6}) gives 
\begin{equation*}
x^{2}=(U_{ms}/U_{m})\left( (V_{ms}^{4}-3V_{ms}^{2}+1)/5\right) .
\end{equation*}%
Clearly, $(U_{ms}/U_{m},(V_{ms}^{4}-3V_{ms}^{2}+1)/5)=1.$ This imples that $%
U_{ms}=U_{m}a^{2}$ and $V_{ms}^{4}-3V_{ms}^{2}+1=5b^{2}$ for some positive
integers $a$ and $b.$ Thus by Theorem \ref{t3.1}, we get $V_{ms}=2,$
implying that $ms=0,$ which is impossible.

Case $2:$ Assume that $5\nmid t.$ Since $5|n,$ it follows that $5|m.$ Then
we can write $m=5^{r}a$ with $5\nmid a,$ $2|a,$ and $r\geq 1.$ It can be
seen by (\ref{2.15}) that $U_{m}=U_{5^{r}a}=5U_{5^{r-1}a}(5a_{1}+1)$ for
some positive integer $a_{1}.$ And thus we conclude that $%
U_{m}=U_{5^{r}a}=5^{r}U_{a}(5a_{1}+1)(5a_{2}+1)...(5a_{r}+1)$ for some
positive integers $a_{i}$ with $1\leq i\leq r.$ Let $%
A=(5a_{1}+1)(5a_{2}+1)...(5a_{r}+1).$ Then we have $U_{m}=5^{r}U_{a}A.$ In a
similar manner, we get $%
U_{n}=U_{5^{r}at}=5^{r}U_{at}(5b_{1}+1)(5b_{2}+1)...(5b_{r}+1)$ for some
positive integers $b_{j}$ with $1\leq j\leq r.$ Let $%
B=(5b_{1}+1)(5b_{2}+1)...(5b_{r}+1).$ It is obvious that $5\nmid B.$ Thus we
have $U_{n}=5^{r}U_{at}B.$ Substituting the new values of $U_{n}$ and $U_{m}$
into $U_{n}=5U_{m}x^{2}$ gives 
\begin{equation}
5^{r}U_{at}B=5.5^{r}U_{a}Ax^{2}.  \label{7}
\end{equation}%
This shows that 
\begin{equation*}
U_{at}B=5U_{a}Ax^{2}.
\end{equation*}%
On the other hand, since $a$ and $at$ are even, it follows from Lemma \ref%
{l2.7} that $U_{at}\equiv (at/2)P(\func{mod}$ $P^{2})$ and $U_{a}\equiv
(a/2)P(\func{mod}$ $P^{2}).$ So (\ref{7}) becomes 
\begin{equation*}
5^{r}(at/2)PB\equiv 5.5^{r}(a/2)PAx^{2}(\func{mod}\text{ }P^{2}).
\end{equation*}%
Rearranging the above congruence gives 
\begin{equation*}
(at/2)B\equiv 5(a/2)Ax^{2}(\func{mod}\text{ }P).
\end{equation*}%
Since $5|P,$ it follows that $5|(at/2)B,$ implying that $5|atB.$ This
contradicts with the fact that $5\nmid a,$ $5\nmid t,$ and $5\nmid B.$ This
completes the proof of Theorem \ref{t3.7}.%
\endproof%


\begin{thebibliography}{99}
\bibitem{KLMN} D. Kalman and R. Mena, \emph{The Fibonacci numbers-exposed,
Mathematics Magazine,} 76 (2003), 167-181.

\bibitem{MSKT} J. B. Muskat, \emph{Generalized Fibonacci and Lucas sequences
and rootfinding methods,} Math. Comp. 61 (1993), 365-372.

\bibitem{RIBENBO} P. Ribenboim,\emph{\ My Numbers, My Friends,}
Springer-Verlag New York, Inc., (2000).

\bibitem{RABINO} S. Rabinowitz, \emph{Algorithmic manipulation of Fibonacci
identities, in Applications of Fibonacci Numbers}, Vol. $6,$ Kluwere, 1996,
pp. 389-408.

\bibitem{LJUN} W. Ljunggrenn, \emph{Zur Theorie der Gleichung }$%
x^{2}+1=Dy^{4},$ Avh. Norsk. Vid. Akad. Uslo, (1942), 1-27.

\bibitem{CO1} J. H. E. John, \emph{On square Fibonacci numbers,} J. London
Math. Soc. (3) 16 (1966), 153-166.

\bibitem{ALF} U. Alfred, \emph{On square Lucas numbers,} Fibonacci Quart. 2
(1964), 11-12.

\bibitem{BURR} S. A. Burr, \emph{On the occurrence of squares in Lucas
sequences,} Amer. Math. Soc. Notices (Abstract 63T-302) 10 (1963), 11-12.

\bibitem{WYLER} O. Wyler, \emph{Solution of problem 5080,} Amer. Math.
Monthly, 71 (1964), 220-222.

\bibitem{CO2} J. H. E. John, \emph{Lucas and Fibonacci numbers and some
Diophantine equations,} Proc. Glasgow Math. Assoc. 7 (1965), 24-28.

\bibitem{CO3} J. H. E. John, \emph{Square Fibonacci numbers, etc.,}
Fibonacci Quart. 2 (1964), 109-113.

\bibitem{CO4} J. H. E. John, \emph{Eight Diophantine equations,} Proc.
London Math. Soc. (3) 16 (1966), 153-166.

\bibitem{CO5} J. H. E. John, \emph{Five Diophantine equations,} Math. Scand.
21 (1967), 61-70.

\bibitem{MC DAN1} P. Ribenboim and W. L. McDaniel, \emph{The square terms in
Lucas sequences,} J. Number Theory 58 (1996), 104-123.

\bibitem{KGW} T. Kagawa and N. Terai, \emph{Squares in Lucas sequences and
some Diophantine equations,} Manuscripta Math. 96 (1998), 195-202.

\bibitem{PETHO} K. Nakamula and A. Petho, \emph{Squares in binary recurrence
sequences,} in: Number Theory, K. Gy\"{o}ry et al. (eds.), de Gruyter,
Berlin, 1998, 409-421.

\bibitem{MC DAN2} P. Ribenboim and W. L. McDaniel, \emph{Squares in Lucas
sequences having an even first parameter,} Collog. Math. 78 (1998), 29-34.

\bibitem{MC DAN3} P. Ribenboim and W. L. McDaniel, \emph{On Lucas sequences
terms of the form }$kx^{2},$ in: Number Theory (Turku, 1999), de Gruyter,
Berlin, 2001, 293-303.

\bibitem{CO6} J. H. E. John, \emph{Squares in some recurrent sequences,}
Pasific J. Math. 41 (1972), 631-646.

\bibitem{KSKN} R. Keskin and Z. Yosma, \emph{On Fibonacci and Lucas numbers
of the form }$cx^{2},$ Journal of Integer Sequences, 14 (2011).

\bibitem{KSKN1} R. Keskin and Z. Yosma, \emph{Some new identities concerning
generalized Fibonacci and Lucas numbers,} Hacettepe Journal of Mathematics
and Statistics, accepted for publication.

\bibitem{NIV} I. Niven, H. S. Zuckerman, and H. L. Montgomary, \emph{An
Introduction to the Theory of Numbers,} John Wiley and Sons, Inc., Toronto,
1991.

\bibitem{BURT} D. M. Burton, \emph{Elementary Number Theory,} Mc Graw--Hill
Group, Inc., New York, 1998.

\bibitem{KSKN2} B. Demirt\"{u}rk and R. Keskin, \emph{Integer solutions of
some Diophantine equations via Fibonacci and Lucas numbers,} Journal of
Integer Sequences, 12 (2009).

\bibitem{ROBBINS} N. Robbins, \emph{On Fibonacci numbers of the form }$%
px^{2},$\emph{\ where }$p$\emph{\ is prime,} Fibonacci Quart. 21 (1983),
266-271.

\bibitem{ROBBINS1} N. Robbins, \emph{Fibonacci numbers of the form }$cx^{2},$%
\emph{\ where }$1\leq c\leq 1000,$ Fibonacci Quart. 28 (1990), 306-315.

\bibitem{ZAFER} Z. \c{S}iar and R. Keskin, \emph{The Square Terms in
Generalized Lucas Sequence,} submitted.

\bibitem{SHOREY} T. N. Shorey and C. L. Stewart, \emph{On the Diophantine
equation }$ax^{2t}bx^{t}y+cy^{2}=1$\emph{\ and pure powers in recurrence
sequences}, Math. Scand. 52 (1983), 24-36.

\bibitem{} R. Keskin and B. Demirt\"{u}rk, \emph{Fibonacci and Lucas
congruences and their applications}, Acta Mathematica Sinica (English
Series), 27 (2011), 725-736.
\end{thebibliography}
\end{document}